\documentclass[11pt]{amsart}
\usepackage{amssymb,amsmath,epsfig,graphics,latexsym,psfrag}
\usepackage[english]{babel}
\usepackage[all]{xy}

\newtheorem{theorem}{Theorem}[section]
\newtheorem{proposition}[theorem]{Proposition}
\newtheorem{corollary}[theorem]{Corollary}
\newtheorem{lemma}[theorem]{Lemma}

\theoremstyle{definition}

\theoremstyle{remark}
\newtheorem{remark}[theorem]{Remark}

\newtheorem{question}{Question}
\numberwithin{equation}{section}
\renewcommand{\t}{ \tilde}
\renewcommand{\b}{ \partial}
\newcommand{\Z}{\bf Z}
\newcommand{\R}{\bf R}

\newcommand{\Hi}{\bf H}

\newcommand{\D}{{\mathcal {D}}}
\renewcommand{\S}{\bf S}
\renewcommand{\l}{\langle}
\renewcommand{\r}{\rangle}
\newcommand{\e}{\varepsilon}
\newcommand{\z}[1]{{\Z}/#1{\Z}}
\renewcommand{\o}{\overline}
\newcommand{\abs}[1]{\lvert#1\rvert}

\newcommand{\co}{\colon\thinspace}
\renewcommand{\epsilon}{\varepsilon}
\renewcommand{\c}{\mathcal}
\usepackage{times}

\begin{document}
\sloppy

\title[]{Finiteness of  mapping degrees  and ${\rm PSL}(2,{\R})$-volume \\on graph manifolds}

\author{Pierre Derbez}
\address{C.M.I,  Technop\^ole Ch\^ateau-Gombert, 39, rue  Joliot Curie, 13453 Marseille Cedex 13, France}
\email{derbez@cmi.univ-mrs.fr}

\author{Shicheng Wang}

\address{Department of Mathematics, Peking University, Beijing, China}

\email{wangsc@math.pku.edu.cn}



\date{\today}
%\maketitle

\begin{abstract} For given closed
orientable $3$-manifolds $M$ and $N$ let $\c{D}(M,N)$ be the set of
mapping degrees from $M$ to $N$. We address the problem: For which
$N$, $\c{D}(M,N)$ is finite for all $M$? The answer is  known for
prime
 3-manifolds
unless the target is a non-trivial graph manifold. We prove that for
each closed non-trivial graph manifold $N$, $\c{D}(M,N)$ is finite
for all graph manifold $M$.

The proof uses a recently developed standard forms of maps between
graph manifolds and the estimation of the $\widetilde{\rm
PSL}(2,{\R})$-volume for certain class of graph manifolds.
\end{abstract}
\maketitle

\vspace{-.5cm} \tableofcontents

\section{Introduction}

Let $M$ and $N$ be two closed oriented $3$-dimensional manifolds.
Let $\c{D}(M,N)$ be the set of degree of maps from $M$ to $N$, that
is
$$\c{D}(M,N)=\{d\in{\Z}\,| f\co M\to N,\, \, \,{\deg}(f)=d\}.$$

According to \cite{CT}, M. Gromov think it is a fundamental problem
in topology to determine the set $\c{D}(M,N)$. Indeed the supremum
of absolute values of degrees in $\c{D}(M,N)$ has been addressed by
J. Milnor and W. Thutston in 1970's \cite{MT}. A basic property of
$\c{D}(M,N)$ is reflected  in the following
\begin{question}\label{gromov} (see  also \cite[Problem A]{Re2} and  \cite[Question 1.3]{W2}):
For which closed orientable $3$-manifold $N$, the set $\c{D}(M,N)$
is finite for all closed orientable $3$-manifolds $M$?
\end{question}

This question can be interpreted as a way to detect some new
rigidity properties of the geometry-topology of a manifold. More
precisely, when $M$ is fixed, then one can expect that if the
geometry-topology of a manifold $N$ is complicated, then the
possible degree of maps $f\co M\to N$ is essentially controled by
the datas of $N$. For geometric $3$-manifolds (i.e. $3$-manifolds
which admits a locally homogeneous complete Riemannian metric) the
answser to this question is summurized in the following
\begin{theorem}[\cite{Th2}, \cite{BG1}, \cite{BG2}, \cite{W2}]\label{known} Let $N$ denote a closed orientable geometric $3$-manifold.

(i) If $N$ supports the hyperbolic or the  $\widetilde{\rm
PSL}(2,{\R})$ geometry, then $\c{D}(M,N)$ is finite for any $M$.

(ii) if $N$ admits one of the six remainder geometry, ${\S}^3$,
${\S}^2\times{\R}$, ${\rm Nil}$, ${\R}^3$, ${\Hi}^2\times{\R}$ or
${\rm Sol}$ then $\c{D}(N,N)$ is infinite.\end{theorem}

To study the set $\c{D}(M,N)$  we need to introduced a special kind
of $3$-manifolds invariants. More precisely, we say that a
non-negative $3$-manifold invariant $\omega$ has \emph{degree
property} or simply \emph{Property D}, if for any map $f\co M\to N$,
 $\omega(M)\geq\abs{{\rm deg}(f)}\omega(N)$. Say $\omega$ has
\emph{covering property} or simply \emph{Property C}, if for any
covering $p\co M\to N$, $\omega(M)=\abs{{\rm deg}(p)}\omega(N)$. The
invariants with Property $D$ are important to study Question 1 due
to the following fact (see Lemma \ref{covering}):

{\bf Fact (*)} {\it If $\omega$ has property $D$ and if $N$ admits a
finite covering $\t{N}$ such that $\omega(\t{N})\not=0$ then the set
$\c{D}(M,N)$ is finite for all $M$}.

When $N$ is hyperbolic, the finiteness of the set  $\c{D}(M,N)$ is
essentially controled by the volume associated to the Riemannian
metric with constant negative sectional curvature which satisfies
\emph{Property D}. When $N$ admits a  $\widetilde{\rm PSL}(2,{\R})$
geometry $\c{D}(M,N)$ is essentially controled by the
$\widetilde{\rm PSL}(2,{\R})$-volume $SV$ introduced in \cite{BG1}
and \cite{BG2}: it satisfies \emph{Property D} and it is non-zero
for  Seifert manifolds supporting $\widetilde{\rm PSL}(2,{\R})$
geometry.

To study Question 1 for more general manifolds, M. Gromov introduced
in \cite{G} the simplicial volume $\|N\|$ of a manifold $N$. This
invariant always satisfies  \emph{Property D}. For example, using
the simplicial volume and the work of Connel Farb in \cite{CF},
Lafont and Schmidt have generalized in \cite{LS} point $(i)$ of
Theorem \ref{known} when the target manifold $N$ is a closed locally
symmetric $n$-manifold with nonpositive sectional curvature.
However, a closed locally symmetric manifold is a special class of
complete locally homogeneous manifolds and thus Question 1 is still
open for non-geometric manifolds with zero simpicial volume.

In this paper we focus on closed $3$-manifolds. Recall that
according to the Perelman geometrization Theorem, $3$-manifolds with
zero Gromov simplicial volume are precisely graph manifolds.   We
call that a  $3$-manifold  covered by either a  torus bundle, or a
Seifert manifold is a \emph{trivial graph manifold}.
 Hence for
prime 3-manifolds, Question \ref{gromov} is reduced to
\begin{question}\label{gromov1}
Suppose $N$ is a non-trivial graph manifold, is $\c{D}(M,N)$  finite
for all closed orientable $3$-manifolds $M$?
\end{question}

The main difficulty to study Question \ref{gromov1} for a
non-trivial graph manifold $N$ is to find a $3$-manifold invariant
satisfying property ${\D}$ which does not vanish of $N$. Based on
Fact (*), it is natural to ask

\begin{question}\label{volume}
Let $N$ be a closed orientable non-trivial graph manifold. Does
$SV(\t{N})\not=0$ for some finite covering $\t{N}$ of $N$?
\end{question}

The $\widetilde{\rm PSL}(2,{\R})$-volume is rather strange and was
very little known. It can be either zero or non-zero for hyperbolic
3-manifolds \cite{BG1}; if it has Property $C$ is still unclear, and
it was not addressed for non-geometric 3-manifolds  since it was
introduced more than 20 years ago.

A main result of this paper is a partial answer of Question
\ref{volume}: we verify that for a family of non-geometric graph
manifolds $N$, they do have finite cover $\t{N}$ with
$SV(\t{N})\not=0$ (Proposition 4.1).  Such a partial answer,
combined with the standard form of nonzero degree maps developed in
\cite{D1}, enable us  to solve Question \ref{gromov1}
 when we restrict on graph manifolds.

\begin{theorem}\label{main}
 For any given closed non-trivial prime graph manifold $N$, $\c{D}(M,N)$ is finite for any graph manifold $M$.
\end{theorem}

\begin{remark}
Some facts related to Theorem \ref{main} are known before:
$\c{D}(N,N)$ is finite for any prime non-trivial graph manifold $N$
(\cite{W1}, see also \cite{D2}). The covering degrees is uniquely
determined by the graph manifolds involved \cite{YW}.
\end{remark}

This paper is organized as follows.

In Section 2 we define the objects which will be used in the paper:
For graph manifolds, we will define their coordinates and gluing
matrices, canonical framings, the standard forms of nonzero degree
maps, the absolute Euler number and the absolute volume.  We also
recall $\widetilde{\rm PSL}(2,{\R})$-volume and its basic
properties.

In section 3 we state and prove some results on coverings of graph
manifolds which  will be used in the paper.

Section 4 is devoted to the proof of Proposition \ref{non-vanish}.
The strategy is to use a finite sequence of coverings  to get a very
"large" and "symmetric" covering space which allows some free action
of a finite cyclic group so that the quotient can be sent onto a
3-manifold supporting $\widetilde{\rm PSL}(2,{\R})$ geometry via a
nonzero degree map.

   In Section 5 we prove Theorem \ref{main}.
 The strategy is to use the standard form of nonzero degree maps
 between graph manifolds to show that one can reduce
 the problem to the case where the target is a graph manifold
 satisfying the hypothesis of Proposition \ref{non-vanish}.

\section{Notations and known results}

From now on all 3-manifolds are irreducible and oriented, and all
graph manifolds are non-trivial.

Suppose $F$ (resp. $P$) is a properly embedded surface (resp. an
embedded 3-manifold) in a 3-manifold $M$.  We use $M\setminus F$
(resp. $M\setminus P$) to denote the resulting manifold obtained by
splitting $M$ along $F$ (resp. removing $\text{int} P$, the interior
of $P$).

\subsection{Coordinated graph manifolds and gluing matrices}
Let $N$ be a graph manifold. Denote by $\c{T}_N$ the family of JSJ
tori of $N$, by $N^{\ast}$ the set $N\setminus
\c{T}_N=\{\Sigma_1,...,\Sigma_n\}$ of the JSJ pieces of $N$, by
$\tau\co\b N^{\ast}\to\b N^{\ast}$ the associated sewing involution
defined in \cite{JS}.

A {\it dual graph}  of $N$, denoted by $\Gamma_N$, is given as
follows: each vertex represents a JSJ piece of $N$; each edge
represents a JSJ torus of $N$; an edge $e$ connects two vertices
$v_1$ and $v_2$ (may be $v_1=v_2$) if and only if the corresponding
JSJ torus is shared by the corresponding JSJ pieces.

Call a dual graph $\Gamma_N$  {\it directed} if each edge of
$\Gamma_N$ is directed, in other words, is endowed with an arrow. Once
$\Gamma_N$ is directed, the sewing involution $\tau$ becomes a well
defined map, still denoted by $\tau\co\b N^{\ast}\to\b N^{\ast}$.

Suppose $N^*$ contains no pieces homeomorphic to $I(K)$, the twisted
$I$-bundle over the Klein bottle.

Let $\Sigma$ be an oriented Seifert manifold which  admits a unique
Seifert fibration, up to isotopy, and $\partial \Sigma\ne
\emptyset$. Denote by $h$  the homotopy class of the regular fiber
of $\Sigma$, by $\c{O}$ the base 2-orbifold of $\Sigma$ and by
$\Sigma^0$ the space obtained from $\Sigma$ after removing the
singular fibers of $\Sigma$. Then $\Sigma^0$ is a ${\S}^1$-bundle
over a surface $\c{O}^0$ obtained from $\c{O}$ after removing the
exceptional points. Then there exists a cross section
$s\co\c{O}^0\to\Sigma^0$. Call $\Sigma$ is {\it coordinated}, if

(1)
such a section $s\co\c{O}^0\to\Sigma^0$ is chosen,

(2) both
$\b\c{O}^0$ and $h$ are oriented so that their product orientation is
matched with the orientation of $\b\Sigma$ induced by that of $\Sigma$.

Once $\Sigma$ is coordinated, then the orientation on $\partial
\c{O}^0$ and the oriented fiber $h$ gives a basis of $H_1(T;{\Z})$
for each component $T$ of $\partial \Sigma$. We also say that
$\Sigma$ is endowed with a $(s,h)$-basis.

 Since $N^{\ast}$ has no
$I(K)$-components then each component $\Sigma_i$ of $N^{\ast}$
admits a unique Seifert fibration, up to isotopy. Moreover each
component $\Sigma_i$ has the orientation induced from $N$. Call $N$
is {\it coordinated}, if each component $\Sigma_i$ of $N^*$ is
coordinated and $\Gamma_N$ is directed.

Once $N$ is coordinated, then each torus $T$ in $\c{T}_N$ is
associated with a unique $2\times 2$-matrix $A_T$ provided by the
gluing map $\tau| \co T_-(s_-, h_-)\to T_+(s_+,h_+)$: where $T_-,
T_+$ are two torus components in $\partial N^*$ provided by $T$,
with basis $(s_-, h_-)$ and $(s_+,h_+)$ respectively, and
$$\tau (s_-, h_-)= (s_+,h_+)A_T.$$

Call $\{A_T, T\in \c{T}\}$ the {\it gluing matrices}.

\subsection{Canonical framings and canonical submanifold}
Let $\Sigma$ denote an orientable Seifert manifold with regular
fiber $h$. A \emph{framing} $\alpha$ of $\Sigma$ is to assign a
simple closed essential curve not homotopic to the regular fiber of
$\Sigma$, for each component $T$ of $\b\Sigma$. Denote by
$\Sigma(\alpha)$ the closed Seifert 3-manifold obtained from
$\Sigma$ after Dehn fillings along the family $\alpha$ and denote by
$\pi_{\Sigma}\co\Sigma\to\Sigma(\alpha)$ the natural quotient map.
Let $p\co\t{\Sigma}\to\Sigma$ be a finite covering. Assume that
$\Sigma$ and $\t{\Sigma}$ are endowed with a framing $\alpha$ and
$\t{\alpha}$. Then we say that $(\t{\Sigma},\t{\alpha})$ covers
$(\Sigma,\alpha)$ if each component of $\t{\alpha}$ is a component
of $p^{-1}(\alpha)$. In this case, the map
$p\co(\t{\Sigma},\t{\alpha})\to (\Sigma,\alpha)$ extends to a map
$\hat{p}\co\t{\Sigma}(\t{\alpha})\to\Sigma(\alpha)$ and the Euler
number of $\Sigma(\alpha)$ is nonzero iff the Euler number of
$\t{\Sigma}(\t{\alpha})$ is nonzero by \cite{LW}. When $N$ contains
no $I(K)$-component in its JSJ-decomposition,  each Seifert piece
$\Sigma$ of $N^{\ast}$ is endowed with a \emph{canonical framing}
$\alpha_{\Sigma}$ given by the regular fiber of the Seifert pieces
of $N^{\ast}$ adjacent to $\Sigma$. Denote by $\hat{\Sigma}$ the
space $\Sigma(\alpha_{\Sigma})$. By minimality of JSJ decomposition,
$\hat{\Sigma}$ admits a unique Seifert fibration extending that of
$\Sigma$.

Call a submanifold $L$ of a graph manifold $N$ is {\it canonical} if
$L$ is a union of some components of $N\setminus {\mathcal T}$,
where ${\mathcal T}$ is subfamily of ${\mathcal T}_N$. Similarly
call $\alpha_L=\{t_U\subset U\}$ where $t_U$ is the regular fiber of
the Seifert piece adjacent to $L$ along the component $U$, when $U$
runs over the components of $\partial L$,
 the canonical framing of $L$, and denote by $\hat L$   the
closed graph 3-manifold obtained from $L$ after Dehn fillings
along the family $\alpha_L$. From the definition we have

\begin{lemma}\label{canonical}
For a given closed graph manifold $M$, there are only finitely many
canonical framed canonical submanifolds  $(L, \alpha_L)$, and thus
only finitely many $\hat L$.
\end{lemma}

\subsection{Standard forms of nonzero degree maps}

We recall here  two  results which are proved in \cite{D1} in a more
general case.  The first result is related to the standard forms  of
nonzero degree maps.
\begin{proposition}(\cite[Lemma 3.4]{D1})\label{haken}
For a given  closed  graph manifold $M$, there is a finite set
${\mathcal H}=\{M_1,...,M_k\}$ of closed graph manifolds satisfying
the following property: for any nonzero degree map $g\co M\to N$
into a closed non-trivial graph manifold $N$ without $I(K)$ piece in
$N^*$, there exists some $M_i$ in ${\mathcal H}$ and a map $f\co
M_i\to N$ such that:

(i) ${\rm deg}(f)={\rm deg}(g)$,

(ii) for each piece $Q$ in $N^{\ast}$, $f^{-1}(Q)$ is a canonical
submanifold of $M$.
\end{proposition}

The following technical "mapping lemma" will be also useful:
\begin{lemma}(\cite[Lemma 4.3]{D1})\label{mapping}
Suppose $f\co M\to N$ is a  map between closed graph manifolds and
$N^*$ contains no $I(K)$ piece. Let $S$ and $S'$ be two components
of $M^{\ast}$ which are adjacent in $M$ along a subfamily ${\mathcal
T}$ of ${\mathcal T}_M$ and satisfy:

(i) $f(S')\subset{\rm int}(\Sigma')$ for some piece $\Sigma'$ of
$N^{\ast}$,

(ii) $f_{\ast}([h_S])\not=1$, where $t_S$ is the regular fiber of
$S$.

Then there exists a piece $\Sigma$ of $N^{\ast}$ and a homotopy of
$f$ supported in a regular neighborhood of $S$ such that
$f(S)\subset{\rm int}(\Sigma)$. Moreover if $f(h_S)$ is not homotopic
to a non-trivial power of the regular fiber of $\Sigma$, then one
can choose $\Sigma=\Sigma'$.
\end{lemma}
\cite[Lemma 4.3]{D1}  was stated for Haken manifolds. Since here we
consider only non-trivial graph manifolds instead of Haken
manifolds, then we can state \cite[Lemma 4.3]{D1} in term of the
JSJ-pieces of $N$ instead of in term of the charactersitic Seifert
pair of $N$.
\subsection{$\widetilde{\rm PSL}(2,{\R})$-volume, absolute volume, and absolute Euler number}

$\widetilde{\rm PSL}(2,{\R})$-volume $SV$ is introduced in
\cite{BG1} and \cite{BG2}. (It is also considered as a special case
of volumes of representations, see \cite{Re1}, \cite{Re2} and
\cite{WZ}). Two basic properties of $SV$ are reflected in the
following

\begin{lemma}\label{SV}
(i) $SV$ has Property $D$ \cite{BG1},

(ii)  If $N$ supports $\widetilde{\rm PSL}(2,{\R})$ geometry, i.e.,
$N$ admits a Seifert fibration with nonzero Euler number $e(N)$ and
whose base $2$-orbifold $\c{O}_N$ has a negative Euler
characteristic, then \cite{BG2}
$$SV(N)=\left|\frac{\chi^2(\c{O}_N)}{e(N)}\right|.$$
\end{lemma}

When $N$ is a closed graph manifold with no $I(K)$ piece in $N^*$,
using the notations introduced in Section 2.2, one can define the
so-called {\it absolute volume $|SV|$} by setting
$$|SV|(N)=\sum_{\Sigma\in N^{\ast}}SV(\hat{\Sigma}).$$
In the same way  one can define the {\it absolute Euler number} of
$N$ by setting
$$|e|(N)=\sum_{\Sigma\in N^{\ast}}|e(\hat{\Sigma})|.$$
In Section 3.3 we will study the relations between $|e|(N)$ and
$|SV|(N)$ (see Lemma \ref{unif}).
\section{Reduction of complexity via coverings}
In this section  we state some results on finite coverings of
surfaces and 3-manifolds which will be used in the proofs of
Proposition 4.1 and Theorem \ref{main}.

\subsection{Two general statements}
The first result says that to prove the finiteness of the set
$\c{D}(M,N)$ one can replace $N$ by a finite covering  of it.

\begin{lemma}\label{covering}
(1) Let $M$ and  $N$ denote two closed oriented manifolds of the
same dimension and let $p\co N'\to N$ be a finite covering of $N$.
If $\c{D}(P,N')$ is finite for any closed manifold $P$ then the set
$\c{D}(M,N)$ is finite.

(2) If all manifolds involved in (1) are graph manifolds, then the
conclusion of (1) is still hold.
\end{lemma}

\begin{proof}
(1) For each nonzero degree map $f\co M\to N$, let $M(f)$ be the connected covering space of $M$ corresponding to the subgroup $f_{\ast}^{-1}(p_{\ast}(\pi_1N'))$ of $\pi_1M$ which we denote by $r\co M(f)\to M$. Let $f'\co M(f)\to N'$ be a lift of $f$, then $p\circ f'=f\circ r$. We claim that the set
$$\c{C}=\{M(f), \ {\rm when}\ f\ {\rm runs\ over\ the\ nonzero\ degree\ maps\ from}\ M\ {\rm to}\ N\}$$
is finite. To see this, first note that the index of  $f_{\ast}^{-1}(p_{\ast}(\pi_1N'))$ in $\pi_1M$ is bounded by the index of $p_{\ast}(\pi_1N')$ in $\pi_1N$. Indeed, the homomorphism $f_{\ast}\co\pi_1M\to\pi_1N$ descends through an injective map $$\bar{f}_{\ast}\co\frac{\pi_1M}{f_{\ast}^{-1}(p_{\ast}(\pi_1N'))}\to\frac{\pi_1N}{p_{\ast}(\pi_1N')}$$ Since $\pi_1M$
contains at most finitely many subgroup of a bounded index, it
follows that $M(f)$ has only finitely many choices which proves that the set $\c{C}$ is finite. By the
construction we have
$${\rm deg}(f)=\frac{{\rm deg}(p)}{{\rm deg}(r)}{\rm deg}(f').$$
By the finiteness of the set $\c{C}$ and assumption on $N'$, the set
of $\{{\rm deg}(f')|f:M'\to N', M'\in \c{C}\}$ is finite.  Clearly
${\rm deg}(r)$ have only finitely many choice, so the lemma is proved.

(2) If $M$ and $N$ are graph manifolds, then all manifolds $M(f)$, $N'$
 in the proof of (1) are graph manifolds. Clearly (2)
follows.
\end{proof}

\begin{lemma}\label{not-shared}
Let $N$ be a closed 3-manifold with non-trivial JSJ-decomposition.
Then there exists a $2$-fold covering $\tilde N$ of $N$ such that
each JSJ-torus of $\tilde N$ is shared by two different pieces of
$\tilde N^*$.
\end{lemma}
\begin{proof} Let $\{T_1,...,T_k\}$ be the union JSJ-tori of $N$ with
each $T_i$ is shared by the same piece of $N^*$. Let $e_1,...,e_k$
be the corresponding edges in $\Gamma_N$. Then $e_1,...,e_k$ are the
edges of $\Gamma_N$ with the  two ends of each $e_i$ being at the
same vertex. Clearly   $ H_1(\Gamma_N;{\Z})= \l e_1\r\oplus...
\oplus\l e_k\r \oplus G$.

\begin{center}
\psfrag{a}{$e_i^2$} \psfrag{b}{$e_i^1$} \psfrag{c}{$e_i$}
\includegraphics[height=1.5in]{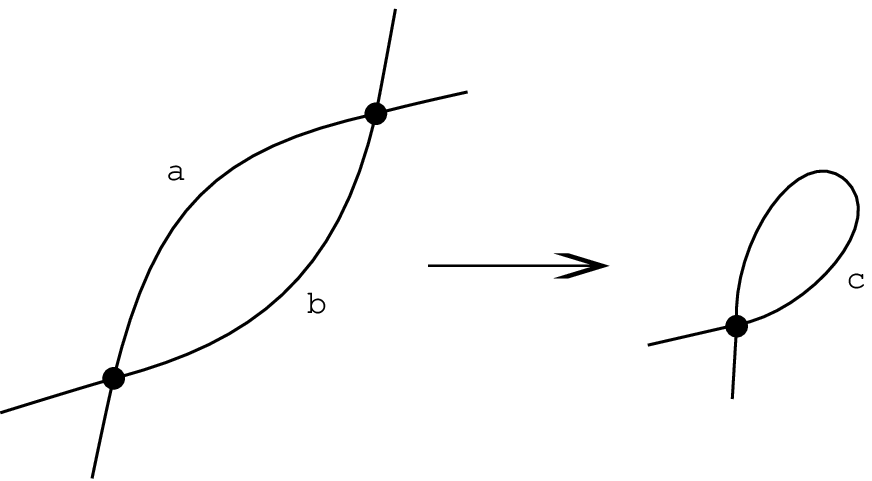}
\centerline{Figure 1}
\end{center}

%\begin{figure}[htb]
%\centerline{\scalebox{.8}{\includegraphics{DWFig1.eps}}} \caption{}
%\end{figure}

Let $r: N \to \Gamma_N$ be the retraction. Consider the following
epimorphism

$$\phi: H_1(N,{\Z}) \stackrel{r_*}\to  H_1(\Gamma_N;{\Z})\stackrel{q}\to\l e_1\r\oplus...
\oplus\l e_k\r \stackrel{\lambda}{\to}\z{2}$$ where $r_*$ is induced
by $r$, $q$ is the projection, and $\lambda$ is defined by
$\lambda([e_i])=\o{1}$ for $i=1,...,k$. Then the double covering
$\t{N}$ of  $N$ corresponding to $\phi$ satisfies the conclusion of
the lemma, since the double covering of $\Gamma_N$ corresponding to
$\lambda\circ q$, which is the dual graph of $\t{N}$, contains no
edge with two ends being at the same vertex. See Figure 1 for the
local picture.
\end{proof}

\subsection{Separable and characteristic coverings}

Let $N$  be a closed graph manifold without  $I(K)$ JSJ-piece.   Let
$\c{T}$ be a union of tori and let $m$ be a positive integer. Call a
covering $p: \t{\c{T}}\to \c{T}$  {\it $m$-characteristic} if for
each component $T$ of $\c{T}$ and for each component $\t{T}$ of
$\t{\c{T}}$ over $T$,  the restriction $p| : \t{T}\to T$  is the
covering map associated to the characteristic subgroup of index $m
\times m$ in $\pi_1T$. Call a finite covering $\t{N}\to N$ of a
graph manifold $N$ $m$-characteristic if its restriction to
$\c{T}_{\t{N}}\to\c{T}_N$ is {\it $m$-characteristic}.

Next we define the \emph{separable coverings}. Let $\Sigma$ be a
component of $N^{\ast}$ with base $2$-orbifold $\c{O}$. Let
$\Sigma^0$, $\c{O}^0$, and $s\co\c{O}^0\to\Sigma^0$ are be given as
in Section 2.1. Let $p\co\t{\Sigma}\to\Sigma$ denote a finite
covering. Recall that $p$ is a fiber preserving map.

Recall that the \emph{vertical degree} of $p$ is the integer $d_v$
such that $p_{\ast}(\t{h})=h^{d_v}$, where $h$ and $\t{h}$ denote
the homotopy class of the regular fiber in $\Sigma$ and
$\t{\Sigma}$, and the \emph{horizontal degree} $d_h$ is the degree
of the induced branched covering $\bar p: \tilde{\c{O}}\to \c{O}$,
where $\t{\c{O}}$ denotes the base of the bundle $\t{\Sigma}$. We
have ${\rm deg}(p)=d_v\times d_h$.

On the other hand, $p$ induces a finite covering
$p|\co\t{\Sigma}^0=p^{-1}(\Sigma^0)\to\Sigma^0$ and a covering
$p|_b\co\t{\c{F}}^0\to\c{O}^0$,  with $\t{\c{F}}^0$ connected. More
precisely, $p|_b$ corresponds to the subgroup
$s_{\ast}^{-1}((p|)_{\ast}(\pi_1\t{\Sigma}^0))$.   Note that $p$ and
$p|$ have the same degree, same  vertical degree and same horizontal
degree. If ${\rm deg}(p|_b)= d_h$, then we say that the covering $p$
is \emph{separable}. The following result provides two classes of
separable coverings which will be used later.

\begin{lemma}\label{separable}
Let  $p\co\t{\Sigma}\to\Sigma$  be an oriented Seifert manifold
finite covering.

(i) If $p$ has fiber degree one, then $p$ is a separable covering.

(ii)  If $\Sigma=F\times{\S}^1$ and
 $p$ is a regular covering corresponding to an epimorphism
$\phi\co\pi_1N=\pi_1F\times{\Z}\to G=G_1\times G_2$ satisfying
$\phi(\pi_1F\times\{1\})=G_1$ and $\phi(\{1\}\times{\Z})=G_2$ then
$p$ is separable.
\end{lemma}
\begin{proof}
Using the same notations as above it is  easy  to see that  the map
$p|_b\co\t{\c{F}}^0\to\c{O}^0$ factors through covering maps
$q\co\t{\c{F}}^0\to\t{\c{O}}^0$ and $\bar p \co\t{\c{O}}^0\to
\c{O}^0$ where $\t{\c{O}}^0$ denote the base of the bundle
$\t{\Sigma}^0$. Then we get
$${\rm deg}(p|_b)=d_h\times{\rm deg}(q)={\rm deg}(p)\times{\rm deg}(q)$$ since $p$ has vertical degree one.
On the other hand, since ${\rm deg}(p|_b)\leq{\rm deg}(p)$ then
${\rm deg}(q)=1$. This proves (i).

If $\Sigma$ is homeomorphic to a product $F\times{\S}^1$ then we have the following commutative diagram
$$\xymatrix{\t{F} \ar[r]^{\t{s}} \ar[d]_{p_b} & \t{\Sigma} \ar[d]^p\\
F \ar[r]^s & \Sigma}
$$ where $\t{F}$ is connected.
Since $\phi(\pi_1F\times\{1\})=G_1$ then $p^{-1}(s(F))$ has $|G_2|$
components and thus ${\rm deg}(p)={\rm deg}(p_b)\times|G_2|$. Since
${\rm deg}(p)=|G_1|\times|G_2|$ then ${\rm deg}(p_b)=|G_1|=d_h$.
This proves (ii).
\end{proof}

\subsection{Lifting of coordinates and gluing matrices}

From now on we assume the graph manifold $N$ is coordinated. Let
$p\co\t{N}\to N$ be a finite covering of graph manifolds. Then
obviously $\Gamma_{\t{N}}$ can be directly in a unique way such that
the induced map $p_\#: \Gamma_{\t{N}}\to \Gamma_N$ preserves the
directions of the edges. Below we also assume that $\Gamma_{\t{N}}$
is directed in such a way.

Let $p\co\t{N}\to N$ be a finite covering of graph manifolds. Call
$p$ is {\it separable} if the restriction $p|
\co\t{\Sigma}\to\Sigma$ on connected Seifert pieces is separable for
all possible $\t{\Sigma}$ and $\Sigma$. Call a coordinate on $\tilde
N$ is a {\it lift} of the coordinate of $N$, if for each possible
covering $p| \co\t{\Sigma}\to\Sigma$ on connected Seifert pieces,
the $(s,h)$-basis of $\t{\Sigma}$ is lifted from the $(s,h)$-basis
of ${\Sigma}$.

\begin{lemma}\label{stable} (i) Let $p\co\t{N}\to N$ be a
separable finite covering of graph manifolds. Then the coordinate of
$N$ can be lifted on $\tilde N$.

(ii) Moreover, if the covering $p$ is characteristic, then for each
component $T$ of $\c{T}_N$ and for each component $\t{T}$ over $T$
we have $A_T=A_{\t{T}}$, where the coordinate of $\tilde N$ is
lifted from $N$.
\end{lemma}
\begin{proof}
To prove (i), one need only to show that for a separable finite
covering $p: \t{\Sigma}\to\Sigma$  of connected Seifert piece, then any
$(s,h)$-basis of $\Sigma$ lifts to a $(s,h)$-basis of $\t{\Sigma}$.

Using the same notation as in the proof of Lemma \ref{separable} we
have
$${\rm deg}(p|_b)=d_h\times{\rm deg}(q)$$
Since we assume  ${\rm deg}(p|_b)= d_h$, then ${\rm deg}(q)=1$ and
thus $\t{s}$ is a cross section. This proves (i).

Once $N$ is coordinated, then each torus $T$ in $\c{T}_N$ is
associated with a unique $2\times 2$-matrix $A_T$ provided by the
gluing map $\tau| \co T_-(s_-, h_-)\to T_+(s_+,h_+)$ such that $\tau
(s_-, h_-)= (s_+,h_+)A_T.$

Similarly with lifted coordinate on $\tilde N$ we have $\tilde \tau|
\co \tilde T_-(\tilde s_-,\tilde h_-)\to\tilde T_+(\tilde s_+,\tilde
h_+)$ and $\tau (\tilde s_-, \tilde h_-)=(\tilde s_+,\tilde
h_+)\tilde A_T.$

Since  the coordinate of $\tilde N$   are lifted from $N$, and
$p$ is $m$-characteristic for some $m$, we have the following
commutative diagram
$$\xymatrix{(\tilde s_-, \tilde h_-) \ar[r]^{\tilde A_T} \ar[d]_{\times m} & (\tilde s_+,\tilde h_+) \ar[d]^{\times m}\\
(s_-, h_-)\ar[r]^{A_T} & (s_+,h_+)}
$$
 Then one verifies directly that $\tilde A_T=A_T$.
This proves (ii).
\end{proof}

\subsection{The absolute volume and the absolute Euler number}
We end this section with a result (see Lemma \ref{unif}) which states the relation between
the absolute volume and the absolute Euler number of a graph manifold. First we begin with a technical result.

\begin{lemma}\label{LW}
Suppose $N$ is a closed graph manifolds without $I(K)$ JSJ-piece.

(i) For any finite covering $\tilde N\to N$, $|e|(\t{N})=0$ if and only
if $|e|(N)=0$.

(ii) There is a finite covering $p:\tilde N\to N$ which is separable
and characteristic, and each Seifert piece of $\tilde N$ is the
product of a surface of genus at least 2 and the circle. Moreover
$\tilde N$ may be chosen so that  $\Gamma_{\tilde N}$ has two
vertices if $\Gamma_N$ has two vertices.
\end{lemma}

\begin{proof}
(i) follows from the definition and \cite[Proposition 2.3]{LW}.

(ii) It has been proved in \cite[Proposition 4.4]{LW}, that there is
a characteristic finite covering $p:\tilde N\to N$  whose each piece
is the product of a surface and the circle. By checking the proof,
it is easy to see that the condition "genus at least 2" can be
satisfied; moreover the restriction $p|: \tilde \Sigma \to \Sigma$
on connected JSJ-pieces is a composition of separable coverings
described in Lemma \ref{separable}, which is still separable.

If moreover $\Gamma_N$ has exactly two vertices $\Sigma_1$ and
$\Sigma_2$, then for $i=1,2$, denote by
$p_i\co\hat{\Sigma}_i\to\Sigma_i$ the $m$-characteristic separable
finite covering such that $\hat{\Sigma}_i$ is the product of a
surface of genus at least 2 and the circle.  There exists a
$1$-characteristic finite covering
$q_i\co\t{\Sigma}_i\to\hat{\Sigma}_i$ such that $\b\t{\Sigma}_1$ and
$\b\t{\Sigma}_2$ have the same number of components. Next one can
glue $\t{\Sigma}_1$ and $\t{\Sigma}_2$ by the lift of the sewing
involution of $N$ to get a characteristic and separable finite
covering $p\co\t{N}\to N$ whose dual graph has two vertices. This
completes the proof of the lemma.
\end{proof}

\begin{lemma}\label{unif}
Let $N$ be a closed graph manifold without $I(K)$ JSJ-pieces.

(i) If $|e|(N)\not=0$ then $N$ admits a finite covering $\t{N}$ with
$|SV|(\t{N})\not=0$.

(ii) If $|e|(N)=0$ then $N$ admits a finite covering $\tilde N$
which can be coordinated such that each gluing matrix is in the form
$\pm\begin{pmatrix}
0&1 \\
1&0
\end{pmatrix}$.
\end{lemma}
\begin{proof}
By Lemma \ref{LW} (ii) let $p:\tilde N\to N$ be a finite covering
which is separable and characteristic and each piece of $\tilde N^*$
is a product $F\times S^1$  with $g(F)\ge 2$.

By Lemma \ref{LW} (i) $|e|(N)\not=0$ implies $|e|(\t{N})\not=0$. By
definition of $|e|$,   $e(\hat{\Sigma})\not=0$ for some $
\Sigma=F\times S^1 \in \tilde N^*$. Since $g(F)\ge 2$,
$SV(\hat{\Sigma})\not=0$, and hence $|SV|(\tilde N)\not=0$ by
definition in Section 2.3. This proves (i).

Denote by $\Sigma_1,...,\Sigma_n$ the components of $N^{\ast}$. For
each $i=1,...,n$, denote by $(\Sigma_i,\alpha_i)$ the  Seifert piece
$\Sigma_i$ of $N^{\ast}$ endowed with its canonical framing. Since
$e(N)=0$ then $e(\Sigma_i(\alpha_i))=e(\hat{\Sigma}_i)=0$ and thus
there exists  a  finite  covering  of $\hat{\Sigma}_i$, with fiber
degree one, homeomorphic to a product. By pulling back this covering
via the quotient map $\pi\co \Sigma_i\to\hat{\Sigma}_i$ we get a
covering $\t{\Sigma}_i$ of $\Sigma_i$ such that the framing
$(\t{\Sigma}_i,\t{\alpha}_i)$ satisfies the following condition:
there exists a properly embedded incompressible surface $F_i$ in
$\t{\Sigma}_i$ such that  $\t{\Sigma}_i\simeq F_i\times{\S}^1$ and
$\b F_i=\t{\alpha}_i$.

Suppose  $T$ is a component of $\b\Sigma_i$ and  $T'$ is a component
of $\b\Sigma_j$ such that $T$ is identified to $T'$ then the sewing
involution $\tau|T\co T\to T'$ lifts to a sewing involution
$\t{\tau}\co\t{T}\to\t{T}'$, where $\t{T}$, resp. $\t{T}'$, denotes
a component of $\b\t{\Sigma}_i$, resp. a component of
$\b\t{\Sigma}_j$, over $T$, resp. $T'$. Indeed by our construction
the induced coverings $\t{T}\to T$ and $\t{T}'\to T'$ correspond
exactly to the subgroup of $\pi_1 T$, resp. of $\pi_1 T'$, generated
by $h$ and $h'$, where $h$ is the fiber of $\Sigma_i$ represented in
$T$ and $h'$ is the fiber of $\Sigma_j$ represented in $T'$, hence
the gluing map lifts by the lifting criterion.

Denote by $\eta_i$ the degree of the covering map
$\t{\Sigma}_i\to\Sigma_i$. Let
$$\eta={\rm l.c.m.}\{\eta_1,...,\eta_n\}$$ For each $i=1,...,n$,
take $t_i=\frac{\eta}{\eta_i}$ copies of $\t{\Sigma}_i$ and glue the
components of
$$\coprod_{i=1}^m\left(\text{$t_i$ copies of }\t{\Sigma}_i\right)$$ together
via lifts of the sewing involution $\tau$ of $N$  to get a separable
finite covering $p\co\t{N}\to N$.
 By coordinating each piece $\tilde \Sigma_i$ of
$\t{N}^*$ with such a section $F_i$ and its regular fiber,
$\tilde N$ is coordinated. Clearly each component of $\b{F}_i$ is
identified with the regular fiber of its adjacent piece and vice
versa. Therefore each gluing matrix should be in the form of
$\begin{pmatrix}
0&\pm1 \\
\pm1&0
\end{pmatrix}$
. Since the determinant should be $-1$, therefore the gluing matrix
is in the form of $\pm\begin{pmatrix}
0&1 \\
1&0
\end{pmatrix}$. This proves  point (ii).
\end{proof}

 \section{$\widetilde{\rm PSL}(2,{\R})$-volume of graph manifolds}

Let $N$ be a closed graph manifold which consists of  two JSJ-pieces
$\Sigma_1$ and $\Sigma_2$ and $n$ JSJ-tori  $\{T_1,...,T_n\}$,
moreover no $\Sigma_i$ is $I(K)$ and each $T_i$ is shared by both
$\Sigma_1$ and $\Sigma_2$. Call such a  manifold
\emph{$n$-multiple edges graph manifold}, whose dual graph $\Gamma_N$ is shown in
Figure 2. We assume $\Gamma_N$ is also directed as in Figure 2. In
this section we use $A_i$ for $A_{T_i}$ for short.

\begin{center}
\psfrag{a}{$T_1$} \psfrag{b}{$T_2$} \psfrag{c}{$T_n$}
\psfrag{d}{$\Sigma_1$} \psfrag{e}{$\Sigma_2$}
\includegraphics[height=2in]{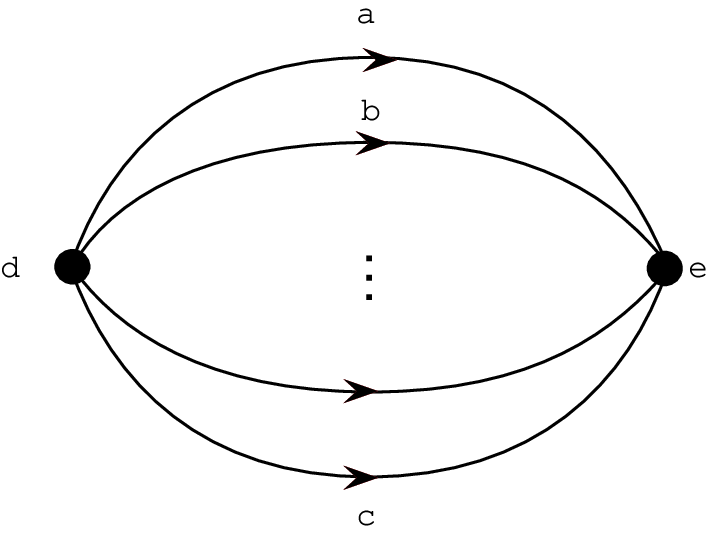}
\centerline{Figure 2: Multiple edges graph manifold}
\end{center}

%\begin{figure}[htb]
%\centerline{\scalebox{.8}{\includegraphics{DWFig2.eps}}}
%\caption{}%Multiple edges graph manifold}
%\end{figure}
\begin{proposition}\label{non-vanish}
Let $N$ be a $n$-multiple edges graph manifold which is coordinated.
Assume that the gluing matrices of $N$ satisfy the condition
$A_{1}=\pm A_{2}=...=\pm A_{n}$. Then $N$ admits a finite covering
space $\t{N}$ such that $SV(\t{N})\not=0$.
 \end{proposition}

\begin{corollary}\label{nonvanish1}
Suppose  $N$ is a closed graph manifold whose dual graph has two
vertices and one edge. Then $|D(M,N)|$ is finite for any 3-manifold
$M$.
\end{corollary}

\begin{proof}
We may suppose that $N$ contains no  $I(K)$ piece.  Otherwise $N$ is
doubly covered by a non-trivial graph manifold which contains no
$I(K)$ piece and whose dual graph still has two vertices and one
edge (since we assume that $N$ is non-trivial graph manifold). In
any case $N$ has a finite cover $\t{N}$ such that $SV(\t{N})\not=0$
by Proposition \ref{non-vanish}. Then by Lemmas \ref{SV} and
\ref{covering}, $|\c{D}(M,N)|$ is finite for any 3-manifold $M$.
\end{proof}

The proof of Proposition \ref{non-vanish} follows from the a
sequence of lemmas below.

\vskip 0.5 truecm

Let $N$ be a $n$-multiple edges coordinated graph manifold. We say
that $N$ sastifies the {\it Property $I$} if

(1) The JSJ-piece $\Sigma_i$ is homeomorphic to a product
$F_i\times{\S}^1$ where $F_i$ is an oriented surface with genus
$\geq 2$, for $i=1,2$;

(2) $A_{1}= A_{2}=...= A_{n}$.

\begin{lemma}\label{I}
Let $N$ be a $n$-multiple edges graph manifold satisfying the
assumption of Proposition \ref{non-vanish}. Then there exist
separable and characteristic finite coverings $p_1\co N_1\to N$ and $p_2\co N_1\to N_2$ such
that $N_2$ satisfies Property $I$.
\end{lemma}
\begin{proof}
By Lemma \ref{LW} (ii) and Lemma \ref{stable}, we may assume that
$N$ is a $n$-multiple edges graph manifold satisfying the assumption
of Proposition \ref{non-vanish}, and moreover $\Sigma_i$ is
homeomorphic to a product $F_i\times{\S}^1$ where $F_i$ is an
oriented surface with genus $\geq 2$.

May assume that $A_1=...=A_k=-A$ and $A_{k+1}=...=A_n=A$, $0<k<n$,
shown as in the left of Figure 3.

Denote by $c_{i,j}$ the loops of $\Gamma_N$ corresponding to the
"composition" $T_i.(-T_j)$, note that here $T_i$ represents an oriented
edge. Then $c_{i,k+1}$ for $i=1,...,k$ and $c_{j,n}$ for
$j=k+1,...,n-1$ form a basis of $H_1(\Gamma_N)$ and we have
$$H_1(\Gamma_N)=(\oplus_{i=1}^k\l c_{i,k+1}\r)\oplus
(\oplus_{j=k+1}^{n-1}\l c_{j,n}\r)$$ Next we define an epimorphism
$$\phi :H_1(N,Z)\stackrel{r_*}\to
H_1(\Gamma_N;{Z})\stackrel{q}{\to}\oplus_{i=1}^k\l
c_{i,k+1}\r\stackrel{\lambda}{\to}\z{2}$$ where $r_*$ is induced by
the retraction $r:N\to \Gamma_N$, $q$ is the projection and
$\lambda$ is defined by $\lambda(c_{i,k+1})=\bar 1$ for any
$i\in\{1,...,k\}$. Denote by $p_1\co N_1\to N$ the $2$-fold covering
corresponding to $\phi$, and by $\mu$ the deck transformation of
this covering.

It is easy to see this covering is  separable and
$1$-characteristic. Moreover with the lifted coordinates of $N$, the
directed graph $\Gamma_{N_1}$ with gluing matrices $\pm A$, as well as the
two lifts $\Sigma_i^1$ and $\Sigma_i^2$ of $\Sigma_i$, $i=1,2$, are
shown in the right of Figure 3.
\begin{center}
\psfrag{a}{$A$} \psfrag{b}{$-A$} \psfrag{c}{$p_1$}
\psfrag{d}{$\Sigma_1$} \psfrag{e}{$\Sigma_2$}
\psfrag{f}{$\Sigma^1_1$} \psfrag{g}{$\Sigma^1_2$}
\psfrag{h}{$\Sigma^2_2$} \psfrag{i}{$\Sigma^2_1$}
\includegraphics[height=2in]{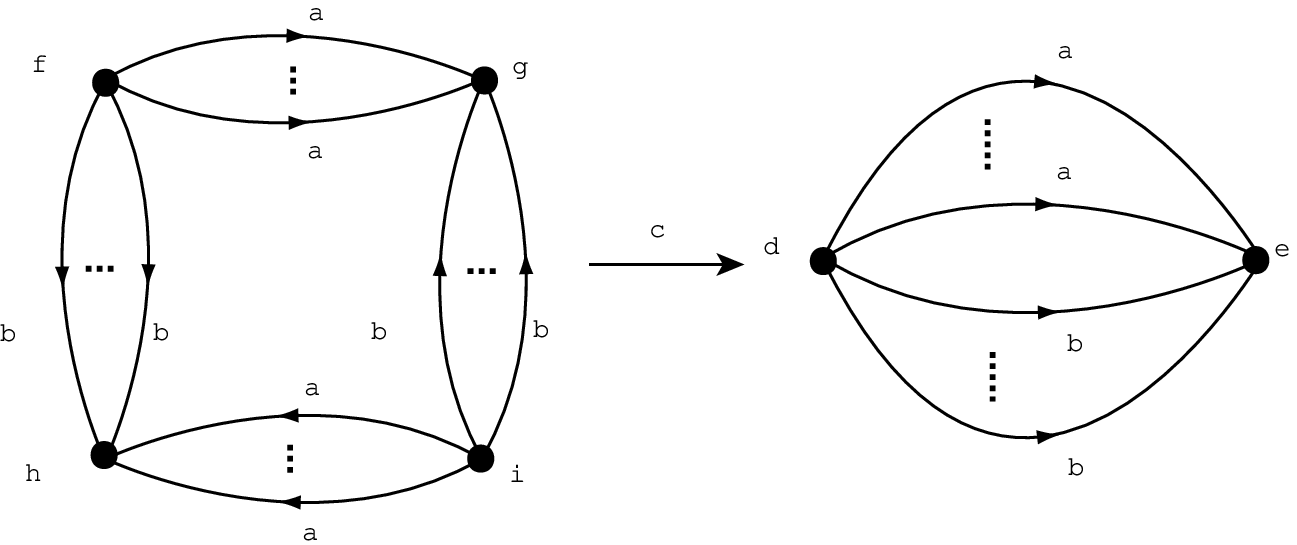}
\centerline{Figure 3}
\end{center}
\vskip 0.5 truecm

%\begin{figure}[htb]
%\centerline{\scalebox{.6}{\includegraphics*[0mm,0mm][210mm,85mm]{DWFig3.eps}}}
%\caption{}
%\end{figure}

Let $\Sigma_i^j= F_i^j\times S^1$. It is not difficult to see that
there is an orientation preserving involution $\eta_i^j$ on
$\Sigma_i^j$ satisfying the following

(1) for $\eta_i^j$ reverses both the orientation of $F_i^j$ and
$S^1$,

(2) for each coordinated component $(T, (s,h))$ of $\partial
\Sigma_i^j$, $\eta_i^j((T, (s,h)))=(T, (-s,-h))$.

Then all those $\eta_i^j$, $i, j=1,2$ match together to get an
involution $\eta$ on $N_1$.

Keep the coordinate of $\Sigma_i^1$ for $i=1,2$, and re-coordinate
$\Sigma_i^2$ for $i=1,2$ by $(T, (-s,-h))$ for each
component of $\partial \Sigma_i^2$ for $i=1,2$, and denoted new
coordinated graph manifold denoted by $N_1'$ ($N_1'$ is $N_1$ if we
forget their coordinates). Then all gluing matrices of $N'_1$ are
$A$.

 Now consider the composition $\eta\circ \mu$, it is easy to see that

(1) $\eta\circ \mu$ is free involution on $N_1'$,

(2) for each JSJ piece $\Sigma$ of $N_1'$, $\eta\circ \mu$ send the
coordinate systems of $\Sigma$ to the the coordinate systems of
$\eta\circ \mu(\Sigma)$.

Now consider the double covering $p_2: N_1=N_1'\to
N_2=N_1'/\eta\circ \mu$. Since the coordinates of $N_1'$ are
invariant under $\eta\circ \mu$, and  all gluing matrices of  $N_1'$
are $A$, we conclude that $N_2$ has Property $I$.
\end{proof}

\begin{lemma}\label{genus}
Let $N$ be a $d$-multiple edges graph manifold satisfying the
property $I$. Then there exists a finite separable
$d$-characteristic covering $p\co N_1\to N$ such that $N_1$ is
$d$-multiple edges graph manifold satisfying the property $I$ and
each JSJ-piece $\Sigma_i^1$ is the product $F^1_i\times{\S}^1$ with
$g(F^1_i)=a_id+b_i\ge 2$ for some positive integers $a_i,b_i$.
    \end{lemma}
    \begin{proof}

Denote by $F_i$ the orbit space of $\Sigma_i$, by $h_i$ its fiber,
and by $c^i_1,...,c^i_{d}$ the components of $\b F_i$ and consider
the homomorphism $$\e_i\co\pi_1\Sigma_i=\pi_1
F_i\times\l[h_i]\r\to\z{d}\times\z{d}$$ defined by
$\e_i(a_l)=(\o{0},\o{0})$ for $l\geq1$, $\e_i(b_j)=(\o{0},\o{0})$
for $j\geq 1$,  $\e_i(c^i_1)=...=\e_i(c^i_{d-1})=(\o{1},\o{0})$ and
$\e_i(h_i)=(\o{0},\o{1})$, where $\pi_1F_i$ has a presentation:
$$\l a_1,b_1,...,a_{g_i},b_{g_i},c^i_1,...,c^i_{d}\ |\ [a_1,b_1]...[a_{g_i},b_{g_i}]c^i_1...c^i_{d}=1\r,$$
where $g_i=g(F_i)$. Since $c^i_1+...+c^i_{d-1}+c^i_{d}=0$ in
$H_1(F_i;{\Z})$ and since $\o{d-1}$ is invertible in $\z{d}$ then
$\e_i(c^i_l)$ is of order $d$ in $\z{d}$ for $l=1,...,d$.  Denote by
$p^1_i\co \Sigma^1_i\to\Sigma_i$ the associated covering, then  the
number of components of $\b\Sigma^1_i$ is $d$  by the construction.
Denote by $p\co N_1\to N$ the $d^2$-fold covering of $N$ obtained by
gluing $\Sigma^1_1$ with $\Sigma^1_2$. This is possible since the
$p^1_i$ induce the $d$-characteristic covering on the boundary for
$i=1,2$. This defines  a finite separable $d$-characteristic
covering by construction. Since $p^1_i$ has horizontal degree $d$
then $\chi(F^1_i)=d\chi(F_i)$, where $F^1_i$ denotes the orbit space
of $\Sigma^1_i$. This implies that
$$2g(F^1_i)+d-2=d(2g(F_i)+d-2)$$ Hence we get
$$g(F^1_i)=d(g(F_i)-1)+\left(\frac{d(d-1)}{2} +1\right).$$ This proves the lemma.
\end{proof}

Call a proper degree one map $p: F'\to F$ between compact surfaces
is a {\it pinch} if there is a  disc $D$ in $\text{int} F$ such that
$p|: p^{-1} (V)\to V$ is a homeomorphism, where
$V=F-\textrm{int}(D)$. We call a proper degree one map $f: F'\times
S^1\to F\times S^1$ a {\it vertical pinch} if $f=p\times \text{id}$,
where $p$ is a pinch.

\begin{lemma}\label{pinch}
Let $N$ be a $d$-multiple edges graph manifold whose gluing matrices
satisfy the condition $I$ and assume that  $g(F_i)=a_id+b_i$ for
some positive integers $a_i,b_i$ and for $i=1,2$. Then $N$ dominates
a $\widetilde{\rm PSL}(2,{\R})$-manifold.
\end{lemma}
\begin{proof}
First note that after perfoming a vertical pinch on
$\Sigma_1=F_1\times\l h_1\r$ and on $\Sigma_2=F_2\times\l h_2\r$ we
may assume that $g_1=g_2=ad+1$ for some $a\in{\Z}_+$.  Note there is
a cyclic $d$-fold  covering ${p}_i': F_i\to  {F_i}'$ with
$g(F_i')=a+1$, $\partial F_i'$ connected, and the restriction of
$p_i'$ is trivial on each component of $\b{F_i}$. This covering is
given by a rotation of angle $2\pi/d$ on $F_i$ whose axis does not
meet $F_i$ (see Figure 4).

\begin{center}
\psfrag{a}{$a_i$ holes} \psfrag{b}{rotation axis} \psfrag{c}{$a_i+1$
holes} \psfrag{1}{$1$} \psfrag{2}{$2$} \psfrag{3}{$3$}
\psfrag{d}{$d$}
\includegraphics[height=5in]{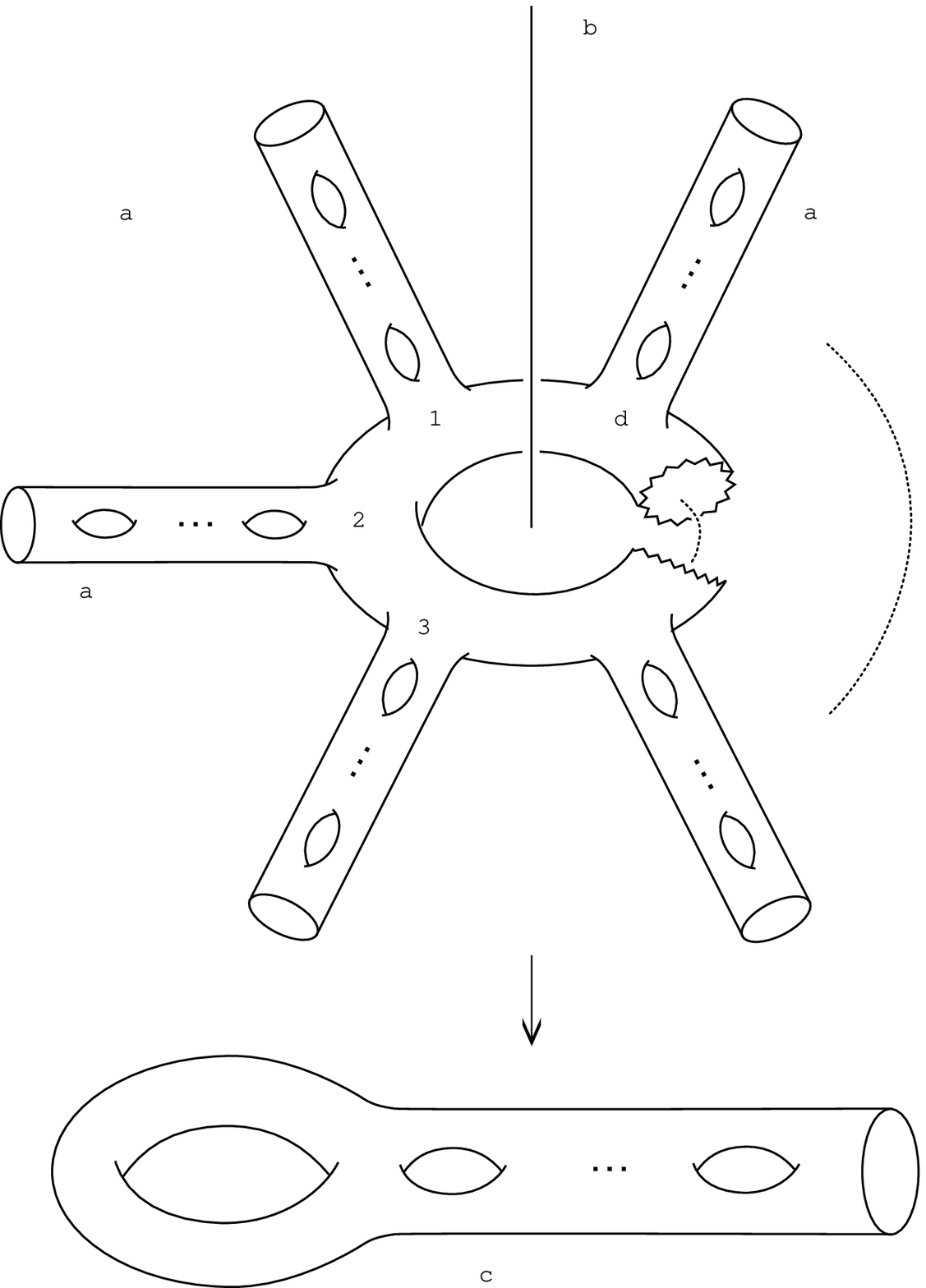}
\centerline{Figure 4: Fixed point free action of $\z{d}$}
\end{center}

%\begin{figure}[htb]
%\centerline{\scalebox{.6}{\includegraphics{DWFig4.eps}}}
%\caption{Fixed point free action of $\z{d}$}
%\end{figure}

The coverings $p_i'\co F_i\to F_i'$ trivially extend to coverings
$p_i'\co\Sigma_i=F_i\times\l h_i\r\to \Sigma_i'=F_i'\times\l h_i'\r$
by setting $p_i'(h_i)=h_i'$ for $i=1,2$. Since all the gluing
matrices of $N$ are $A$ by Property $I$,  then the coverings
$p_i'\co\Sigma_i\to\Sigma_i'$ extend to a covering $p'\co N\to N'$,
where the graph manifold $N'$ consists of the Seifert pieces
$\Sigma_1'$ and $\Sigma_2'$ and the gluing matrix $A$ under obvious
basis.

We fix some notations. For $i=1,2$, denote by $\b{F}_i'=s_i$ and
${\tau}'\co\b{\Sigma}_1'={s}_1\times{h}_1'\to\b{\Sigma}_2'={s}_2\times{h}_2'$
the induced sewing map satisfying
${\tau}'({s}_1,{h}_1')=({s}_2,{h}_2')A$, where $A=\begin{pmatrix}
a&b \\
c&d\end{pmatrix}\in{\rm SL}(2,{\Z})$ with $ad-bc=-1$. Moreover
$b\not=0$ by the basic properties of JSJ-decomposition.

Note that ${\tau}_{\ast}'(s_1)=a {s}_2 +c{h}_2'$ and
${\tau}_{\ast}'^{-1}(s_2)=-d{s}_1+c{h}'_1$. If $ac\not=0$, then
first pinch ${\Sigma}_1'={F}_1'\times {h}_1'$ into a solid torus
$V_1={ D}^2\times{h}_1'$ by killing ${s}_1$. This pinch provides a
degree one map $\pi\co{N}'\to\hat {\Sigma}_2'$ where $\hat
{\Sigma}_2'$ is the closed $3$-manifold obtained from ${\Sigma}_2'$
by Dehn filling along the curve $a {s}_2 +c{h}_2'$. Since $ac\not=0$
then $\Sigma_2'$ is a $\widetilde{\rm PSL}(2,{\R})$-manifold. If
$dc\not=0$ similarly one can perform the same construction with
${\Sigma}_1'$. This proves the lemma when $ac\not=0$ or $dc\not=0$.

Let us assume now that $ac=dc=0$. Then either  $c=0$ or  $a=d=0$.
Since $ad-bc=-1$, then either $A= \pm \begin{pmatrix}
1&b \\
0&-1
\end{pmatrix}$ with $b\not=0$ or $A=\pm\begin{pmatrix}
0&1 \\
1&0
\end{pmatrix}$.

Keeping the coordinate on $\Sigma_1'$ and re-coordinating
$\Sigma'_2$ by $(-s_2,-h_2')$ if needed, we may assume that
$A=\begin{pmatrix}
1&b \\
0&-1
\end{pmatrix}$ or $A=  \begin{pmatrix}
0&1 \\
1&0
\end{pmatrix}$.

Suppose first $A=\begin{pmatrix}
1&b \\
0&-1
\end{pmatrix}$. Denote  ${F}_1'\simeq{F}_2'$ by $F$. Let  denote by
$$\pi\co{F}_1'\times{h}_1'\coprod {F}_2'\times
{h}_2'\to \Sigma=F\times h$$  the trivial 2-fold covering map, where
$\pi({h}_i')=h$ and $\pi(s_i)=s$. Denote by $\rho\co\Sigma\to\hat{\Sigma}$ the
quotient map associated with the Dehn filling on $\b\Sigma$ along
the curve $\frac{b}{(2,b)}s-\frac{2}{(2,b)}h$, where $(2,b)$ denotes
the greatest common divisor of $2$ and $b$. Note that, since
$b\not=0$, then $\hat{\Sigma}$ is a $\widetilde{\rm
PSL}(2,{\R})$-manifold.

One can verify routinely that in the $\pi_1$ level the relations
provided by gluing $\Sigma'_1$ and $\Sigma_2'$ via $\tau'$ are sent
to the relation provided by Dehn filling on $\Sigma$ via
$\frac{b}{(2,b)}s-\frac{2}{(2,b)}h$ under $\rho$, hence the map
$\rho\circ\pi\co {\Sigma}_1'\coprod {\Sigma}_2'\to{\hat \Sigma}$
factors through $\frac{{\Sigma}_1'\coprod
{\Sigma}_2'}{{\tau}'}\simeq{N}'$ which is sent into $\hat {\Sigma}$
by a degree 2 map, since the sewing involution ${\tau}'$ is
orientation reversing so that ${N}'$ inherits compatible
orientations from the pieces ${\Sigma}_1'$ and ${\Sigma}_2'$.

In the case of $A=\begin{pmatrix}
0&1 \\
1&0
\end{pmatrix}$, we can perform the same construction as above, just
replace the filling curve $\frac{b}{(2,b)}s-\frac{2}{(2,b)}h$ by the
curve  $s-h$. This proves Lemma \ref{pinch}.
\end{proof}

By Lemmas \ref{I}, \ref{genus} and \ref{pinch} and their proofs, we
have the following diagram:

  $$\xymatrix{
   N_1 \ar[d]_{p_{1}} \ar[dr]_{p_{2}} & N_3 \ar[d]_{p_{3}} \ar[dr]_{p_{4}}   \\
   N & N_2 & N_4 }$$
where $p_1$ and $p_2$ are coverings provided by Lemma \ref{I}, $p_3$
is the coverings provided by Lemma \ref{genus}, and the non-zero
degree map $p_4$ is provided by Lemma \ref{pinch}, where
$SV(N_4)\not=0$. Since $SV$ has property $D$, $SV(N_3)\not=0$.

Consider the covering $\t{N}$ corresponding to the finite index
subgroup $p_{2\ast}(\pi_1 N_1)\cap p_{3\ast}(\pi_1 N_3)$ in $\pi_1
N_2$. Then $\t{N}$ covers both $N_1$ (and thus $N$) and $N_3$, and
$SV(\tilde N)\not=0$. Then the proof of Proposition
\ref{non-vanish} is complete.

\section{Proof of Theorem \ref{main}}

\subsection{Simplifications}

Let $N$ be a closed non-trivial graph manifold. We are going to show that
$|D(M,N)|$ is finite for any given graph manifold $M$.

(1) First we simplify $N$: By Lemma \ref{not-shared} and Lemma
\ref{LW}, there is a finite covering $\tilde N$ of $N$ satifying the
condition (*): each JSJ piece of $\tilde N$ is a product of an oriented
surface with genus  $\ge 2$ and the circle, and each JSJ torus is
shared by two different JSJ pieces.

By Lemma \ref{covering}, if $|D(M,N)|$ is not finite for some graph
manifold $M$, then $|D(P,\tilde N)|$ is not finite for some graph
manifold $P$. So we may assume $N$ already satisfies the condition
(*).

(2) Then we simplify $M$: For given $M$, let ${\mathcal H}=\{M_1, ...
, M_k\}$ be the finite set of graph manifolds provided by
Proposition \ref{haken}. By Proposition \ref{haken} (i), if
$|D(M,N)|$ is not finite, then $|D(M_i,N)|$ is not finite for some
$M_i\in {\mathcal H}$. So may assume that (**) $M=M_i\in {\mathcal
H}$ for some $i\in\{1,...,k\}$.

\subsection{Proof of Theorem \ref{main} when $|e|(N)\ne 0$.}

Suppose $|e|(N)\not=0$. By Lemma \ref{unif} and (*) in 5.1, we may
assume that $|SV|(N)\not=0$. Then there exists a Seifert piece $Q$
of $N^{\ast}$ such that $SV(\hat{Q})\not=0$. By (**) in 5.1 and
Propostion \ref{haken} (ii), we may assume that $L_Q(f)=f^{-1}(Q)$ is
a canonical submanifold of $M$. Below we denote $L_Q(f)$ as $L_Q$
for short.

\vskip 0.5 true cm

\begin{lemma}\label{different}  $L_Q$ can be chosen so
that any component $T$ of $\b L_Q$ is shared by two distinct Seifert
pieces of $M$: one in $L_Q$ and another in $M\setminus L_Q$.
\end{lemma}

\begin{proof} Indeed if not, then there exists two distinct components $T$
and $T'$ of $\b L_Q$ which are identified by the sewing involution $\tau_M$ of $M$ and
such that $T$ and $T'$ are sent by $f$ into the same component of
$\b Q$. Denote by $\bar{L}_Q$ the canonical submanifold of $M$
obtained by identifying $T$ and $T'$ via $\tau_M$. Since each
component of $\b Q$ is shared by two distinct Seifert pieces of $N$
by the assumption in 5.1,  $f$ induces a proper map
$\bar{f}\co\bar{L}_Q\to Q$. After finitely many such operations, we
reach a new $L_Q$ satisfying the requirement of Lemma
\ref{different}
\end{proof}

Below we assume that $L_Q$ satisfies the requirement of Lemma
\ref{different}. Now we choose $L_Q$ to be maximal in the sense that
for any Seifert piece $S$ in $M\setminus L_Q$ adjacent to $L_Q$, $S$
is not able to be added into $L_Q$ by homotopy on $f$. Then
$f(S)\subset B_{S}$, where $B_S$ is a Seifert piece of $N$, distinct from $Q$ and adjacent
to $Q$.

Since $L_Q$ is maximal, by Lemma \ref{mapping}, we deduce that for
any Seifert piece $S$ adjacent to $L_Q$ along a component of $\b
L_Q$, $f|S\co S\to B_S$ is fiber preserving. Hence the proper  map
$f|L_Q\co L_Q\to Q$ preserves the canonical framings, and it induces
a map $\hat{f}\co\hat{L}_Q\to\hat{Q}$ between the closed manifolds
obtained after Dehn filling along the canonical framings.
 By Lemma \ref{SV}
 we have
$$SV(\hat{L}_Q)\geq\abs{{\rm deg}(\hat{f})}SV(\hat{Q}).$$

Since ${\rm deg}(f)={\rm deg}(f|L_Q)={\rm deg}(\hat{f})$, we get
$$|{\rm deg}(f)|\leq\frac{SV(\hat{L}_Q)}{SV(\hat{Q})}.$$

Therefore

$$|{\rm deg}(f)|\leq\max
\left\{\frac{SV(\hat{L})}{SV(\hat{Q})}| Q\in N^*, SV(\hat Q)\ne 0; L
\,\text{is canonical in}\ M\right\}$$

By Lemma \ref{canonical} there are only finitely many $\hat{Q}$ and   only finitely
many $\hat L$.
So the right side of the above inequality is finite. This completes
the proof of Theorem \ref{main} when $|e|(N)\not=0$.

\subsection{Proof of Theorem \ref{main} when $|e|(N)= 0$.}
By Lemma \ref{covering} and Lemma \ref{unif} we can assume each
gluing matrix of $N$ is equal to $\pm\begin{pmatrix}
0&1 \\
1&0
\end{pmatrix}$.

Choose two distinct  adjacent Seifert pieces $S_1$ and $S_2$ in $N$,
denote by $\c{T}=\b S_1\cap\b S_2$  and by $Q$ the connected graph
manifold $S_1\cup_{\c{T}}S_2$ (such Seifert pieces exist by 5.1). By
Proposition \ref{haken}, we may assume that $f^{-1}(Q)=L_Q$ is a
canonical submanifold of $M$.

Since each JSJ-torus of $N$ is shared by two different JSJ-pieces,
by the same arguments as in 5.2, we may assume that each component of
$\b L_Q$ is shared by two distinct Seifert pieces of $M$ one in
$L_Q$ and another in $M\setminus L_Q$. Furthermore we can arrange
$L_Q$ to be maximal in the sense of 5.2, then by Lemma
\ref{mapping}, we deduce that  any Seifert piece $S'$ of $M$
adjacent to $L_Q$ is sent by $f$ to a Seifert piece $B'$ adjacent to
$Q$ such that $f|S'\co S'\to B'$ is fiber preserving.

As in 5.2,  it follows that the proper map $f|L_Q\co L_Q\to Q$ induces a map
 $\hat{f}\co\hat{L}_Q\to\hat{Q}$ between closed graph manifolds
obtained by Dehn filling along the canonical framings. Moreover, as in 5.2 we have ${\rm deg}(f)={\rm deg}(f|L_Q)={\rm deg}(\hat{f})$ and thus
$$|\c{D}(M,N)|\leq\max\left\{|\c{D}(\hat{L},\hat{Q})|  L
\,\text{is canonical in}\ M\right\}$$ Note that $\hat Q=
\hat{S}_1\cup_{\c{T}} \hat{S}_2$, where $\hat{S}_i$ is obtained by
Dehn filling along the canonical framings on $\b S_i \setminus
\c{T}$, $i=1,2$.  It follows that $\hat{Q}$ satisfies the hypothesis
of Proposition \ref{non-vanish}. Then $\hat{Q}$ has a finite
covering $\tilde Q$ with $SV(\tilde Q)\ne 0$ by Proposition
\ref{non-vanish}. Hence  by Lemma \ref{covering}, the set
$|\c{D}(\hat{L},\hat{Q})|$ is finite for any  $\hat{L}$. Since  by
Lemma \ref{canonical} there are only finitely many $\hat{L}$, this
completes the proof of Theorem \ref{main} when $|e|(N)\not=0$. Hence
 Theorem \ref{main} is proved.

\end{document}